\documentclass[12pt]{amsart}
\usepackage{amsmath}
\usepackage{amsfonts}
\usepackage{amssymb}
\usepackage[latin1]{inputenc}
\usepackage{latexsym}
\usepackage{amscd}
\usepackage{graphicx}
\usepackage{color}

\input amssym.def
\newcommand {\sobre} {\longrightarrow}
\newcommand {\tende} {\rightarrow}

\newcommand {\C} {\mathbb{C}}

\newcommand {\N} {\mathbb{N}}

\newtheorem{definition}{Definition}[section]

\newtheorem{proposition}[definition]{Proposition}
\newtheorem{theorem}[definition]{Theorem}
\newtheorem{corollary}[definition]{Corollary}
\newtheorem{lemma}[definition]{Lemma}
\newtheorem{example}[definition]{Example}

\newtheorem{remark}[definition]{Remark}

\begin{document}
\title{The algebra of bounded type holomorphic functions on the ball}
\author{Daniel Carando}
\address{Departamento de Matem\'{a}tica,
	Facultad de Cs. Exactas y Naturales, Universidad de Buenos Aires and IMAS-UBA-CONICET, Argentina.}
\author{Santiago Muro}
\address{Departamento de Matem\'{a}tica
Facultad de Cs. Exactas y Naturales, Universidad de Buenos Aires, Argentina, and CIFASIS-CONICET}
\author{Daniela M. Vieira}
\address{Departamento de Matem\'{a}tica, Instituto de Matem\'atica e Estat\'istica, Universidade de São Paulo, São Paulo, Brasil}
\email{dcarando@dm.uba.ar} \email{muro@cifacis-conicet.gov.ar} \email{danim@ime.usp.br}
%\subjclass[2000]{Primary 47H60, 47L20. Secondary 46G25, 46M05}

\thanks{The first author was  supported by CONICET-PIP 11220130100329CO, ANPCyT PICT 2015-2299 and UBACyT
20020130100474BA. The second author was  supported by CONICET-PIP 11220130100329CO, ANPCyT PICT 2015-2224 and
	UBACyT 20020130300052BA. The third author was supported by FAPESP-Brazil, Proc. 2014/07373-0}

%\thanks{This work was partially supported by CONICET, FAPESP-Brazil, Proc. 2014/07373-0}

\begin{abstract}
	We study the spectrum $M_b(U)$ of the algebra of bounded type holomorphic functions on a complete Reinhardt domain in a symmetrically regular Banach space $E$ as an analytic manifold over the bidual of the space. In the case that $U$ is the unit ball of $\ell_p$, $1<p<\infty$, we prove that each connected component of $M_b(B_{\ell_p})$ naturally identifies with a ball of a certain radius. We also provide estimates for this radius and in many natural cases we have the precise value. As a consequence, we obtain that for connected components different from that of evaluations, these radii are strictly smaller than one, and can be arbitrarily small.  We also show that for other Banach sequence spaces, connected components do not necessarily identify with balls.
\end{abstract}

\maketitle

\section{Introduction}

The study of the spectrum of the algebra of bounded type analytic functions on a Banach space $E$ was initiated by the seminal article of Aron, Cole and Gamelin \cite{arcoga}. Their main motivation was its relation with the algebra $\mathcal H^\infty(B_E)$ of bounded holomorphic functions on the unit ball. As in the one or finite dimensional case, there is a natural projection defined on the spectrum $\mathcal M$ of  $\mathcal H^\infty$, which in the infinite dimensional case, has range contained in the closed unit ball of the bidual $\overline{B_{E''}}$.

%Moreover, in the finite dimensional case, this projection is one-to-one over the open unit ball and the only (great) difficulties come from the study of the homomorphisms that are projected to the boundary of the unit ball. The famous Corona Problem, proved to be true by Carleson in the one dimensional case, and which remains open for the $n$-dimensional case ($n\ge2$), asks whether the unit ball (or its preimage under the projection) is dense in $\mathcal M$ endowed with the $w^*$ (or pointwise convergence) topology.
%In contrast, in the infinite dimensional case, the spectrum $\mathcal M$ may be highly non-trivial also over the interior of the unit ball and the Corona Problem is also an open problem for the homomorphisms which are mapped to the open unit ball  $B_{E''}$.

The results proved in \cite{arcoga} imply that the interior part of the spectrum $\mathcal M$ (i.e. the subset of homomorphisms which lie in the fibers of the interior points of the ball)  naturally identifies with the spectrum $M_b(B_{E})$ of the algebra of bounded type holomorphic functions on the unit ball of the Banach space $E$.

In \cite{argagama}, the authors continued the study of the spectrum of the algebra of bounded type analytic functions. They showed that for symmetrically regular Banach spaces, the spectrum $M_b(U)$ of the algebra $H_b(U)$ of bounded type holomorphic functions on an open set $U\subset E$ may be endowed with an analytic structure as an infinite dimensional Banach manifold modeled over the bidual $E''$ of $E$. This was applied, for example, to characterize the envelope of holomorphy of $U$ in \cite{carmur, muro}. The analytic structure of $M_b(X)$ for $X$ a Riemann domain over a symetrically  regular Banach space was studied in \cite{dinven}.

\medskip
In this article, we study  the spectrum of the algebra of bounded type analytic functions on the unit ball   of $E$ (or on a complete Reinhardt domain) from this point of view. More precisely, we aim to give an accurate description of $M_b(U)$ as analytic manifold. We show that whenever $U$ is a complete Reinhardt domain in a reflexive space with 1-unconditional basis, each connected component of $M_b(U)$ is (identified with) a complete Reinhardt set, which is not necessarily a multiple of $U$. We also prove that, when $U$ is the unit ball of $\ell_p$, the connected components are identified with balls in the following sense (see definitions below): they are all of the form
\begin{equation}\label{eq-identif}
S= \{\varphi^z\ : \ \|z\|<r\},
\end{equation}
for some $\varphi$ in the fiber of 0 and some $0<r\le 1$. Moreover, with the exception of the  component formed by evaluations, the radius $r$ is strictly smaller that 1.
Also, there are connected components with arbitrary small radius. To show these facts, we give estimates of the radius of each connected component and, for the components of most natural homomorphisms, we give their exact value. This altogether provides a thorough description of $M_b(B_{\ell_p})$, which in turn gives information on the spectrum of $\mathcal H^{\infty}(B_{ \ell_p})$ by \cite{arcoga}.

The fact that connected components are identified with balls as in \eqref{eq-identif} is a particular (isometric) property of $\ell_p$: we  exhibit an example of a Banach space $E$ with 1-unconditional basis for which  the connected components of $M_b(B_E)$ are not balls. The example is actually a Banach space isomorphic to $\ell_2$.

\medskip
We refer to \cite{din, muj} for general theory on complex analysis in Banach spaces, and to \cite{cagama, carandosurvey,din,gamelinlecture, daniindag} for  background on the space of holomorphic functions of bounded type and its spectrum.

%\textcolor[rgb]{0.00,1.00,0.00}{
%CITAR?: LIBRO DINEEN, SURVEY CARANDO etal,  ZAGORODNYUK?, VIEIRA, DINEEN-VENKOVA, LECTURES DE ARON Y DE GAMELIN?.... Y ALGO MAS?
%}

\section{The spectrum of bounded type functions on complete Reinhardt domains}

Let $E$ be a complex Banach space. We denote by $E'$ its dual, and by $B_E$ its open unit ball. Sometimes, when the underlying Banach space is clear, we use $B_r(x)$ to denote the open ball of radius $r$ centered at $x$ and write $B_r$ when the ball is centered at the origin.

For an open subset $U\subset E$, a $U$-\emph{bounded} set is a bounded set $A\subset U$ whose distance to the boundary of $U$, denoted by $d_U(A)$, is positive. {A family $(U_n)_{n\in\N}$ of subsets of $U$ is a \emph{fundamental family} of $U$-bounded sets if each $U_n$ is $U$-bounded, and if every $U$-bounded set is contained in some $U_n$. Every open set $U$ admits a fundamental family of $U$-bounded sets, for instance $$U_n=\{x\in U\,:\,\|x\|\le n,\, d_U(x)\ge \frac 1 n\},$$ for every $n\in\N$.
A holomorphic function on $U$ which is bounded on $U$-bounded sets is called of \emph{bounded type} on $U$.
The algebra of all bounded type holomorphic functions on $U$ is denoted by $H_b(U)$ and it is a Fr\'echet algebra when it is endowed with the topology of uniform convergence on $U$-bounded sets. The \emph{spectrum} of $H_b(U)$, i.e. the set of non-zero continuous complex valued homomorphisms on $H_b(U)$, is denoted by $M_b(U)$.
For each homomorphism $\varphi\in M_b(U)$, there exists a $U$-bounded subset $A$ such that
\begin{equation}\label{varphi<A}
|\varphi(f) | \leq \Vert f \Vert_{A}, \quad\textrm{ for every }f\in H_b(U),
\end{equation}
where $\|f\|_{A}$ is the supremum of $|f|$ over the set $A$. We will write $\varphi\prec A$ when \eqref{varphi<A} holds.

There is a natural projection $\pi:M_b(U)\to E''$, defined by $\pi
(\varphi)=\varphi |_{E'} \in E''$, $\varphi \in M_b(U)$. We thus have the following commutative diagram:
\vspace{0.5cm}
\begin{center}
	\begin{picture}(130,90)
	\put(0,85){$U$} \put(110,85){$M_b(U)$} \put(120,15){$E''$}
	\put(15,88){\vector(1,0){95}} \qbezier(15,88)(12,88)(12,91) \qbezier(12,91)(12,94)(15,94)
	\put(123,81){\vector(0,-1){55}} \put(10,80){\vector(2,-1){110}}
	\qbezier(10,80)(8.4,80.8)(9.2,82.4) \qbezier(9.2,82.4)(10,84)(11.6,83.2) \put(62,90){$\delta$}
	\put(126,50){$\pi$}\put(55,45){$j_E$}
	\end{picture}
\end{center}
where $\delta$ is the point evaluation mapping and $j_E:E\to E''$ is the natural inclusion.

A Banach space $E$ is \emph{symmetrically regular} if every continuous symmetric linear mapping $T:E \to E'$ is weakly
compact (an operator $T:E\to E'$ is symmetric if $Tx_1(x_2) = Tx_2(x_1)$ for all $x_1,x_2 \in E$). Every reflexive Banach space is symmetrically regular.
In \cite{argagama}, for $E$ a symmetrically regular Banach space and $U\subset E$ an open subset, a topology is defined on $M_b(U)$ so that the mapping $\pi$ above is a local homeomorphism that makes $(M_b(U),\pi)$ a Riemann domain over $E''$.

Let us briefly describe this topology (see \cite{argagama} for details). Recall that any  holomorphic function $f$ of bounded type on $E$ may be extended to a function $AB(f)\in H_b(E'')$ through the Aron-Berner extension \cite{arober}.
Given $f\in H_b(U)$ and $z\in E''$, the function
$$
x\mapsto AB\big(\frac{d^nf(x)}{n!}\big)(z),
$$
is a bounded type holomorphic function on $U$.
For
$\varphi\in M_b(U)$, we denote by $d_U(\varphi)$ the supremum of $d_U(A)$ over the $U$-bounded sets $A$ satisfying  $\varphi \prec A$. If $r<d_U(\varphi)$, it is possible to define, for each $z\in E''$ with  $\|z\|<r$, the homomorphism $\varphi^z$ given by
\begin{equation}\label{fi supra z}
\varphi^z(f)=\sum_{n=0}^\infty
	\varphi\Big(AB\big(\frac{d^nf(\cdot)}{n!}\big)(z)\Big).
\end{equation}
When $E$ is symmetrically regular, the sets $\{\varphi^z:\; \|z\|<r\}$, with $\varphi\in M_b(U)$ and $r<d_U(\varphi)$, form a basis of a Hausdorff topology for $M_b(U)$, and each set $\{\varphi^z:\; \|z\|<r\}$ is homeomorphic to the ball $\pi(\varphi)+r  B_{E''}$ via the projection $\pi$. This endows $M_b(U)$ with an analytic structure over $E''$.

\begin{definition}
	Let $U$ be an open subset of a symmetrically regular Banach space. The connected component of a homomorphism $\varphi\in M_b(U)$ is called the \emph{sheet} of $\varphi$ in $M_b(U)$ and is denoted by $S_U(\varphi)$.
\end{definition}

In the case of bounded type entire functions (i.e., $U=E$), the description of the connected components of $M_b(E)$ is simpler than for a general open set $U$, as pointed out in \cite{argagama} and \cite[Section 6.3]{din}.  Given $z\in E''$, the function $x\mapsto\tau_zf(x):=AB(f)(z+j_Ex)$ is an entire function of bounded type on $E$.  Thus, given $\varphi\in M_b(E)$ and $z\in E''$, the homomorphism $\varphi^z$ can be equivalently constructed  as $$\varphi^z(f):=\varphi(\tau_zf).$$ The \emph{sheet} of $\varphi$ is exactly $$S_E(\varphi):=\{\varphi^z\,:\,z\in E''\}.$$ Since $\pi(\varphi^z)=\pi(\varphi)+z$,  $\pi$ is a homeomorphism between $S_E(\varphi)$ and $E''$.
%Thus, these sheets are exactly the connected components of $M_b(E)$.

%%\textcolor[rgb]{0.00,1.00,0.00}{EXPLICAR MEJOR ESTO. DECIR QUE $\pi$ ES LA PROYECCION DESDE EL ESPECTRO DE LAS ENTERAS AUNQUE SE LO APLIQUEMOS A ELEMENTOS DE $M_b(U)$. PROPONER LA NOTACION $\varphi\in M_B(U)\cap \pi^{-1}(0)$ O ALGUNA SIMILAR}\textcolor{red}{ahi lo explique muy resumidamente! no entiendo bien lo de la notacion, vos decis que notemos esto $\varphi\in M_b(U)\cap \pi^{-1}(0)$ de otra manera? porque es muy largo?}
\begin{remark}\rm \label{U balanceado}
If $U\subset E$ is a balanced open set  (or more generally, if $U$ is such that entire functions of bounded type are dense in $H_b(U)$), the spectrum $M_b(U)$ is naturally embedded in $M_b(E)$. Indeed, given $\varphi\in M_b(U)$ we can naturally associate a unique character on $H_b(E)$ which is just the restriction to the bounded type entire functions: $\varphi_{|_{H_b(E)}}$. When the context is clear we will denote this restriction by $\varphi_|$.
The natural projection defined on $M_b(U)$ is just the restriction of the projection defined on $M_b(E)$, and we will denote both as $\pi$.
\end{remark}

Suppose that $U$ is balanced. The embedding of $(M_b(U),\pi)$ into $(M_b(E),\pi)$ is continuous (with their topologies as Riemann domains), so each connected component of $M_b(U)$ is embedded into a connected component of $M_b(E)$ (which is homeomorphic to $E''$).
Therefore, restricted to each sheet, the projection $\pi|_{S_U(\varphi)}$ is a homeomorphism onto some open set of $E''$. Our main goal is to describe the connected components $S_U(\varphi)$, and a natural way to do this is to understand the image  $\pi|_{S_U(\varphi)}$.

Under the same assumptions,
given $\varphi\in M_b(U)$ and $\psi\in S_U(\varphi)$
there exists %$\varphi\in M_b(E)$ belonging to $\pi^{-1}(0)$
  $z\in E''$ such that $\psi_|= (\varphi_|)^z$ and then $(\varphi_|)^z$ belongs to $M_b(U)$ (that is, it can be extended to $H_b(U)$).
Thus, to describe what the connected components of $M_b(U)$ look like, it will be useful to determine for which $z\in E''$ the homomorphism $(\varphi_|)^z$ belongs to $M_b(U)$ (which means, again, that $(\varphi_|)^z$ can be extended to $H_b(U)$).

\smallskip

The following lemma from \cite{aclm} will be useful for our results, in particular for Lemma \ref{lema1}.
\begin{lemma}\label{aclm} \cite[Lemma 1.7]{aclm}. Let $E$ be a Banach space with Schauder basis $(e_k)_{k\in\N}$, and denote by $(e_k')_{k\in\N}$ its dual basic sequence. Let $z\in E''$ and $\varphi\in\pi^{-1}(z)$. Then for $f\in H_b(E)$ and $N\in\N$: $$\varphi(f)=\varphi\Big(x\mapsto f\big(\sum_{k=1}^N z(e_k')e_k+\sum_{k=N+1}^{\infty}e_k'(x)e_k \big) \Big).$$\end{lemma}

\begin{lemma}\label{lema1}Let $E$ be a Banach space with Schauder basis $(e_k)_{k\in\N}$,  and let $\varphi\in M_b(E)\cap \pi^{-1}(0)$. For each $N\in\N$, the following assertions hold.
\begin{enumerate}
\item[\textbf{(1)}] For $z\in E''$ and  $f\in H_b(E)$,  $$\varphi^z(f)=\varphi(x\mapsto AB(f)(z_1,\dots,z_N,x_{N+1}+z_{N+1},x_{N+2}+z_{N+2},\dots)).$$
\item[\textbf{(2)}] If $\varphi\prec A$, then $\varphi\prec A^{(N)}$, where $$A^{(N)}=\{(0,\dots,0,x_{N+1},x_{N+2},\dots)\,:\,x=(x_j)\in A\}.$$
\end{enumerate}
\end{lemma}

\textbf{Proof:} If $\varphi\in\pi^{-1}(0)$, then $z=0$ in Lemma \ref{aclm}. Then $$\varphi(f)=\varphi(x\mapsto f\big(\sum_{k=N+1}^{\infty}e_k'(x)e_k\big))=\varphi(x\mapsto f(0,\dots,0,x_{N+1},x_{N+2},\dots)).$$

\textbf{(1)}  If $f\in H_b(E)$, then $\varphi^z(f)=\varphi(x\mapsto AB(f)(x+z))$. If we denote $g(x)=AB(f)(x+z)$, then it follows from Lemma \ref{aclm} that $$\varphi^z(f)=\varphi(x\mapsto g(0,\dots,0,x_{N+1},x_{N+2},\dots)=AB(f)(z_1,\dots,z_N,x_{N+1}+z_{N+1},x_{N+2}+z_{N+2},\dots)).$$

\textbf{(2)} Since $\varphi\prec A$, we have  $$|\varphi(f)|=|\varphi(x\mapsto f(0,\dots,0,x_{N+1},x_{N+2},\dots))|\le \sup_{x\in A}|f(0,\dots,0,x_{N+1},x_{N+2},\dots)|=\sup_{A^{(N)}}|f|.\  \Box$$

%Recall that given $\varphi\in M_b(U)$, the \emph{sheet} of $\varphi$, $S(\varphi)$, is the connected component of $M_b(U)$ that contains $\varphi$.

We recall that a subset $U$ of a Banach space with unconditional basis $(e_k)_{k\in\N}$ is \emph{complete Reinhardt} if $\sum_{k=1}^{\infty}\lambda_kx_k e_k\in U$, whenever $\sum_{k=1}^{\infty}x_k e_k\in U$ and $|\lambda_k|\leq 1$ for all $k$. Proposition \ref{sheet1} states that if $U$ is a complete Reinhardt domain in a Banach space with 1-unconditional basis, then each sheet in the spectrum is also a complete Reinhardt domain. First we need the following lemma, which is probably known.

\begin{lemma}\label{reinhardt fundamental}
Let $E$ be a Banach space with unconditional basis and let $U\subset E$ be a complete Reinhardt open set. Then $U$ admits a fundamental system of $U$-bounded sets formed by complete Reinhardt sets.
\end{lemma}
\textbf{Proof:}
Any Banach space with unconditional basis can be renormed so that $\|\lambda\cdot x\| \le\|x\|$
whenever $\|\lambda\|_\infty\le 1$. Assuming that $E$ has such a norm, let us show that the
sets $U_n=\{x\in U\,:\,\|x\|\le n,\, d_U(x)\ge \frac 1 n\}$ are complete Reinhardt.
Note that it suffices to prove that if $B_\delta(x)\subset U$ and $\|\lambda\|_\infty\le1$, then $B_\delta(\lambda\cdot x)\subset U$.

Let $y$ be a point in $B_\delta(\lambda\cdot x)$ and define a vector $z \in E$ by specifying its
coordinates as follows:
$$
z_j=\left\{
\begin{array}{ll}
y_j, & \textrm{ if } x_j=0,\\
\frac{x_j}{|x_j|}\max(|x_j|,|y_j|) & \textrm{ otherwise. }
\end{array}
\right.
$$
	If the index $j$ is such that $|x_j | < |y_j |$, then $|z_j - x_j | = |y_j | - |x_j | \le |y_j | - |\lambda_j x_j | \le
	|y_j -\lambda_j x_j |$ by the triangle inequality. And if $j$ is such that $|x_j |\ge |y_j |$, then $|z_j -x_j | = 0 \le
	|y_j - \lambda_j x_j |$. Thus $|z_j - x_j | \le |y_j - \lambda_j x_j |$ for every index $j$, so $\|z - x\|\le\|y -\lambda\cdot x\| < \delta$.
	In other words, $z \in B_\delta (x)$, so $z\in U$. Since $|z_j |\ge |y_j |$ for every $j$, and $U$ is a complete
	Reinhardt set, it follows that $y\in U$. But $y$ is an arbitrary point of $B_\delta (\lambda\cdot x)$, so we conclude that $B_\delta (\lambda\cdot x) \subset
	U$.\hfill$\Box$\medskip

If we only look at the subset of homomorphisms that project to $E$, then the above topology restricted to $M_b(U)\cap \pi^{-1}(E)$ is well defined, even though $E$ is not symmetrically regular. Thus, for an arbitrary Banach space $E$,  $(M_b(U)\cap\pi^{-1}(E),\pi|_{\pi^{-1}(E)})$ is a Riemann domain over $E$ (see \cite{carmur}).

\begin{proposition}\label{sheet1}
Let $E$ be Banach space with 1-unconditional basis $(e_k)_{k\in\N}$ and let $U\subset E$ be a complete Reinhardt open subset.
Then, in each sheet of $M_b(U)$ there is a character $\varphi\in M_b(U)\cap\pi^{-1}(0)$ such that the set
$$
\{w\in E \,:\,(\varphi_{|})^w\textrm{ extends to } M_b(U)\}
$$
is a complete Reinhardt subset of $E$.
\end{proposition}

\textbf{Proof:}
Recall that since $H_b(E)$ is dense in $H_b(U)$, we have that $M_b(U)$ is embedded in $M_b(E)$. Then, given $\psi\in M_b(U)\cap \pi^{-1}(E)$
there exists $\varphi\in M_b(E) \cap \pi^{-1}(0)$ and $z\in E$ such that $\psi_|=\varphi^z$.
%$S_U(\psi)\cap\pi^{-1}(0)$ is embedded in $S_E(\varphi)\cap\pi^{-1}(0)$.
We must show that for every scalar sequence $\lambda$ with   $\|\lambda\|_\infty\le1$, the vector $w=\lambda\cdot z$ satisfies that $\varphi^w$ extends to $M_b(U)$ whenever $\varphi^z$ extends to $M_b(U)$.
  Note that since $\varphi^w$ belongs to $M_b(E)$, it suffices to show that $\varphi^w\prec A$ for some $U$-bounded set $A$.

Let us start by assuming that $z=\sum_{j=1}^Nz_je_j$. %and $\lambda_j=0$ for every $j>N$.
If $f\in H_b(E)$, it follows by Lemma \ref{lema1} that
$$
\varphi^w(f)=\varphi(x\mapsto f(\lambda_1 z_1,\dots,\lambda_N z_N,x_{N+1},x_{N+2},\dots)).
$$
Let us consider the entire function of bounded type,
$$f_{\lambda}(x)=f(\lambda_1 x_1,\dots,\lambda_N x_N,x_{N+1},x_{N+2},\dots),$$
then, applying again  Lemma \ref{lema1},
$$\varphi^z(f_{\lambda})=\varphi(x\mapsto f_\lambda(z_1,\dots,z_N,x_{N+1},x_{N+2},\dots))=\varphi(x\mapsto f(\lambda_1 z_1,\dots,\lambda_N z_N,x_{N+1},x_{N+2},\dots))=\varphi^w(f).$$
By the previous lemma we may take a complete Reinhardt $U$-bounded set, $A$, such that $\varphi^z\prec A$.
Then,
$$|\varphi^w(f)|=|\varphi^z(f_{\lambda})|\leq\sup_A|f_{\lambda}|\leq \sup_A|f|.$$
Therefore $\varphi^w\in M_b(U)$ and $\varphi^w\prec A$.

Take now an arbitrary $z\in E$ for which  $\varphi^z$ belongs to $M_b(U)$ with $\varphi^z\prec A$. Let us denote by $\pi_N$ the projection onto the span of $\{e_1,\dots,e_N\}$ and choose $0<\delta<\frac{d_U(A)}{3}$. We can take $N$ such that $\|\pi_N(z)-z\|<\delta<\frac{d_U(A)}{3}$. Now,  proceeding as in \cite[page 550]{argagama}, we have $\varphi^{\pi_N(z)}\prec A_\delta:=A+B_{\delta}$. By the first part of the proof, for $\|\lambda\|_\infty\le1$ we have $\varphi^{\lambda\cdot\pi_N(z)}\prec A_\delta$. Since $d_U(A_\delta)>2\delta$ and $\|\lambda\cdot\pi_N(z)-\lambda\cdot z\|<\delta$, we have $\varphi^{\lambda\cdot z}\prec A_{2\delta}$. Finally, since $\delta$ is arbitrary small,  we conclude that $\varphi^{\lambda\cdot z}\prec A$.
\hfill$\Box$\medskip

%
%
%\begin{corollary}\label{sheet1cor}
%Let $E$ be a reflexive Banach space with 1-unconditional basis and let $U\subset E$ be a complete Reinhardt open subset.
%Then (the projection of) each sheet is a complete Reinhardt set. That is, for each sheet $S$ of $M_b(U)$ there is a character $\varphi\in M_b(U)\cap\pi^{-1}(0)$ such that $S=S_U(\varphi)$ and
%$$
%\pi(S)=\{w\in E\,:\,\varphi^w\in M_b(U)\}
%$$
%is a complete Reinhardt subset of $E$.
%\end{corollary}
%
%
%\begin{remark}
%Another way to state the previous corollary is the following: for each sheet $S$ of $M_b(U)$ there exists $\varphi\in M_b(U)\cap\pi^{-1}(0)$ and a complete Reinhardt domain $V\subset E$ such that
%$$
%S= \{\varphi^z\in M_b(U)\ : \ z\in V\}.
%$$
%\end{remark}

\smallskip

If the Banach space $E$ is reflexive (which obviously implies that $E$ is symmetrically regular), the above result tells us that the sheets of $M_b(U)$ are complete Reinhardt domains.

\begin{corollary}\label{sheet1cor}
	Let $E$ be a reflexive Banach space with 1-unconditional basis and let $U\subset E$ be a complete Reinhardt open subset.
	Then for each sheet $S$ of $M_b(U)$ there exist a character $\varphi\in M_b(U)\cap\pi^{-1}(0)$ and a complete Reinhardt domain $V\subset E$ such that
	$$
	S= \{(\varphi_|)^z\in M_b(U)\ : \ z\in V\}.
	$$
\end{corollary}

\section{The spectrum of bounded type functions on $B_{\ell_p}$}

We now focus in the case where $U$ is the unit ball of $\ell_p$. The following theorem shows that each sheet is also a ball centered at zero. We will see later in Theorem  \ref{sheet3} that the radius of each sheet other than the sheet of evaluations, is strictly smaller than 1.

\begin{theorem}\label{sheet2} Let $E=\ell_p$, $1<p<\infty$, and let $U=B_{\ell_p}$. Then all sheets are balls centered at 0, that is, in each sheet there is some $\varphi\in M_b(U)\cap\pi^{-1}(0)$, and
$$
\pi(S_U(\varphi))=\{w\in E\,:\,(\varphi_|)^w\in M_b(U)\}=rB_{\ell_p},
$$
for some $0<r\le1$.
%. If $\varphi^z\in S(\varphi)$ and if $w\in\ell_p$ is such that $\|w\|<\|z\|$, then $\varphi^w\in S(\varphi)$.
\end{theorem}
%\begin{theorem}\label{sheet2} Let $E=\ell_p$, $1<p<\infty$, and let $U=B_{\ell_p}$. Then all sheets are balls centered at 0, that is, in each sheet $S$ there is some $\varphi\in M_b(U)\cap\pi^{-1}(0)$, and some $0<r\le1$ such that
%	$$
%	S=S_U(\varphi)= \{\varphi^z\ : \ \|z\|<r\}.
%	$$
%\end{theorem}
\textbf{Proof:} By Corollary \ref{sheet1cor} we know that each sheet intersects $\pi^{-1}(0)$. So take $\varphi\in M_b(U)\cap\pi^{-1}(0)$ and suppose that $(\varphi_|)^z$ belongs to $M_b(U)$ for some $z\in E$. The theorem will be proved if we show that $(\varphi_|)^w\in M_b(U)$ whenever $\|w\|<\|z\|$.

If $w=(w_j)_{j\in\N}$ and $z=(z_j)_{j\in\N}$ are such that $\|w\|<\|z\|$, then there exists $N_1\in\N$ such that $\|\sum_{j=1}^{N}w_j e_j\|< \|\sum_{j=1}^{N}z_j e_j\|$ for every $N\ge N_1$. On the other hand, since $(\varphi_|)^z\in S_U(\varphi)$,  there exists $\delta>0$ such that $(\varphi_|)^{z+y}\in S_U(\varphi)$, for all $\|y\|<\delta$. So let us take $N\geq N_1$ such that
$$
\sum_{j=N+1}^{\infty}|z_j|^p<\Big(\frac{\delta}{3}\Big)^p \textrm { and } \sum_{j=N+1}^{\infty}|w_j|^p<\Big(\frac{\delta}{3}\Big)^p.
$$

Then, if $v=(\Pi_N(z),(I-\Pi_N)(w))$, where $\Pi_N:\ell_p\sobre\ell_p$ denotes the canonical projection, we have that  $(\varphi_|)^{v}$ also belongs to $S_U(\varphi)$.
Note that $\|w\|<\|v\|$ and that $(I-\Pi_N)(w)=(I-\Pi_N)(v)$.

To show that $(\varphi_|)^w\in S_U(\varphi)$, we will construct some auxiliary bounded linear transformations, as follows.
First, take $\gamma:\C^N\sobre\C$ such that $\|\gamma\|=\|(v_1,\dots,v_N)\|^{-1} $ and $\gamma(v_1,\dots,v_N)=1$.
Next, we define $S_N:\C^{N}\sobre\C^{N}$ by $$S_N(x)=\gamma(x)(w_1,\dots,w_N),$$ which clearly satisfies $\|S_N\|\leq 1$ and $S_N(v_1,\dots,v_N)=(w_1,\dots,w_N)$.   Finally, let $T_N:\ell_p\sobre\ell_p$ be given by $T_N(x)=(S_N(\Pi_N(x)),(I-\Pi_N)(x))$. In other words, $$T_N(x)=(S_N(x_1,\dots,x_N),x_{N+1},x_{N+2},\dots),\quad\textrm{for }x \in\ell_p.$$ Note that $T_N(v)=w$ and, since
\begin{eqnarray*}
\|T_N(x)\|^p&=&\|S_N(\Pi_N(x))\|^p+\|(I-\Pi_N)(x)\|^p  \leq \|S_N\|^p\|\Pi_N(x)\|^p+\|(I-\Pi_N)(x)\|^p \\ &\leq& \|\Pi_N(x)\|^p+\|(I-\Pi_N)(x)\|^p=\|x\|^p,\end{eqnarray*} we also have $\|T_N\|\leq 1$.

%Define $T_N:\ell_p\sobre\ell_p$ by $T_N(x)=(S_N(\Pi_N(x)),(I-\Pi_N)(x))$, that is,
%$$
%T_N(x)=(S_N(x_1,\cdots,x_N),x_{N+1},x_{N+2},\cdots),\textrm{ for all }x \in\ell_p,
%$$
%where $S_N:\C^{N}\sobre\C^{N}$ is such that $\|S_N\|\leq 1$ and $S_N(z_1,\cdots,z_N)=(w_1,\cdots,w_N)$. The mapping $S_N$ can be constructed, for example, in the following way:  take $\gamma:\C^N\sobre\C$ such that $\|\gamma\|=\|(z_1,\cdots,z_n)\|^{-1} $ and $\gamma(z_1,\cdots,z_n)=1$, then $S_N(x)=\gamma(x)(w_1,\cdots,w_N)$ has the required conditions.
%
%Note that $T_N(v)=w$. Moreover, $\|T_N\|\leq 1$, indeed,
% $\|T_N(x)\|^p=\|S_N(\Pi_N(x))\|^p+\|(I-\Pi_N)(x)\|^p\leq\|S_N\|^p\|\Pi_N(x)\|^p+\|(I-\Pi_N)(x)\|^p\leq \|\Pi_N(x)\|^p+\|(I-\Pi_N)(x)\|^p=\|x\|^p$.

If $f\in H_b(E)$, then it follows from Lemma \ref{lema1} that $$(\varphi_|)^v(f)=\varphi(x\mapsto f(\Pi_N(v),(I-\Pi_N)(x+v))$$ and that $$(\varphi_|)^w(f)=(\varphi_|)^{T_N(v)}(f)=\varphi(x\mapsto f(\Pi_N(T_N(v)),(I-\Pi_N)(x+T_N(v))).$$
Since  $\Pi_N(T_N(v))=S_N(\Pi_N(v))$ and $(I-\Pi_N)(T_N(v))=(I-\Pi_N)(v)$, we have
\begin{equation}\label{eq1}
(\varphi_|)^{w}(f)=(\varphi_|)^{T_N(v)}(f)=\varphi(x\mapsto f(S_N(\Pi_N(v)),(I-\Pi_N)(x+v)).
\end{equation}

On the other hand, for $f\in H_b(B_{\ell_p})$, consider $f_N=f\circ T_N|_{B_{\ell_p}}\in H_b(B_{\ell_p})$. Then we have
$$
\begin{array}{rl}
	f_N(\Pi_N(v),(I-\Pi_N)(x+v)) & =f\circ T_N\big(\Pi_N(v),(I-\Pi_N)(x+v)\big) \\
	&  =f(T_N(v_1,\dots,v_N,x_{N+1}+w_{N+1},x_{N+2}+w_{N+2},\dots))\\
	&=f(S_N(\Pi_N(v)),(I-\Pi_N)(x+v)).\end{array}
$$
Hence, \begin{eqnarray*}(\varphi_|)^v(f_N)&=&\varphi(x\mapsto f_N(\Pi_N(v),(I-\Pi_N)(x+v))=\varphi(x\mapsto f(S_N(\Pi_N(v)),(I-\Pi_N)(x+v)) \\ &=&(\varphi_|)^w(f).\end{eqnarray*}

If $A$ is a $U$-bounded ball such that $\varphi^{v}\prec A$, then, using again that $\|T_N\|\leq 1$, we conclude   that
$$
|(\varphi_|)^w(f)|%=|\varphi^{T_N({z})}(f)|
=|(\varphi_|)^{{v}}(f_N)|\leq\sup_A|f_N|=\sup_{T_N(A)}|f|\leq\sup_A|f|,$$
which shows that $(\varphi_|)^w\in S_U(\varphi)$.
  \hfill$\Box$\medskip.

%We want to use Theorem \ref{sheet2} to show that when $E=\ell_p$, $1<p<\infty$, $U=B_{\ell_p}$ and $\varphi\in M_b(U)$ is such that $\varphi\in\pi^{-1}(0)$, then $\pi(S(\varphi))$ is an open ball centered at zero. For that, we need the following lemma.
%
%\begin{lema}\label{lemabola} Let $A$ be a nonempty open and bounded subset of a Banach space $E$ with the following property: if $z\in A$ and $\|w\|<\|z\|$ then $w\in A$. Then $A$ is an open ball centered at zero. \end{lema}
%
%\textbf{Proof:} Let $r=\sup_{x\in A}\|x\|$. Then it is not difficult to show that $B(0,r)\subset A \subset B[0,r]$, and then we conclude that $A=B(0,r)$. \hfill$\Box$\medskip.
%
%\begin{corollary} Let $E=\ell_p$, $1<p<\infty$, $U=B_{\ell_p}$ and let $\varphi\in M_b(U)$ be such that $\varphi\in\pi^{-1}(0)$. Then $\pi(S(\varphi))$ is an open ball centered at zero.\end{corollary}
%
%\textbf{Proof:} We observe that $\pi(S(\varphi))$ is bounded, see, for instance, \cite[Propostion 18]{cagama}. Then just apply Lemma \ref{lemabola}. \hfill$\Box$\medskip.
%
%\bigskip

A natural question at this point is whether each sheet on $M_b(B_E)$  is necessarily a ball centered at zero, for more general Banach spaces. The next example shows that this is not always true.
\begin{example}\rm
 Let $E=\langle e_0\rangle\oplus_\infty\ell_2$. Take $\varphi\in M_b(B_{E})$ to be any limit point of the sequence $(\delta_{e_n/\sqrt{2}})_n$.
 By Proposition \ref{sheet1} we know that the projection of the sheet of $\varphi$
 $$
 \pi(S_{B_E}(\varphi))=\{x\in E\ :\ (\varphi_|)^x\in M_b(B_E)\},
 $$
 is a complete Reinhardt open set. Let us show that $\pi(S_{B_E}(\varphi))$ is not a ball centered at 0. For this we will see that $(\varphi_|)^{se_0}\in M_b(B_E)$ for every $|s|<1$ but that $(\varphi_|)^{te_1}\notin M_b(B_E)$ for every $|t|>1/\sqrt 2$.

 For the first assertion, just note that the set $({se_0+e_n/\sqrt{2}})_n$ is $B_E$-bounded and clearly $(\varphi_|)^{se_0}\prec({se_0+e_n/\sqrt{2}})_n$, thus $(\varphi_|)^{se_0}\in M_b(B_E)$. For the second assertion, define the function $f(x)=\sum_{k\ge1}x_k^2$. Then $f\in H_b(E)$ and for every $m\in\mathbb N$, its $m^{th}$-power satisfies $\|f^m\|_{B_E}=1$. On the other hand, since  $(\varphi_|)^{te_1}\in M_b(E)$, we know that for each $m\in\mathbb N$,  $(\varphi_|)^{te_1}(f^m)$ is a limit point of $(\delta_{te_1+e_n/\sqrt{2}}f^m)_n$. Finally, since $f(te_1+e_n/\sqrt{2})^m=(t^2+\frac12)^m\to \infty$ as $m\to\infty$, we conclude that $(\varphi_|)^{te_1}$ cannot be extended to $H_b(B_E)$.
\end{example}

Now that we know that each sheet of $M_b(B_{\ell_p})$ is a ball centered at zero, we would like to estimate its radius. Let us first recall some terminology from \cite{arcoga} that will be used in the next theorem. For $\varphi\in M_b(B_E)$ and $m\ge0$ we associate  $\varphi_m\in \mathcal P(^mE)'$, as $\varphi_m:=\varphi|_{\mathcal P(^mE)}$. Recall also that $R(\varphi)$, the radius of $\varphi$, is defined as the infimum of all $r>0$ such that $\varphi\prec rB_E$. In \cite{arcoga} it is shown that $$R(\varphi)=\limsup_{m\in\N}\|\varphi_m\|^{\frac{1}{m}}=\sup_{m\in\N}\|\varphi_m\|^{\frac{1}{m}}.$$
It should be mentioned that the definition of the radius and the above result were given for $\varphi\in M_b(E)$, but it is easily checked that the same works for $\varphi\in M_b(B_E)$.

%We will use the following Polarization Formula, which appeared in \cite[Lemma 1.7]{CarDimMur09} for the case $k=1$. We include a short proof for sake of completeness.
%\begin{lemma}\label{polarization}
%Let $\rho$ be a primitive $n$th-root of 1 (i.e. $\rho^n=1$  and  $\rho^j\neq 1$ for $j<n$), $P$ an $n$-homogeneous polynomial and $1\le k<n$. Then
%$$\check{P}(z^k,x^{n-k})=\dfrac{1}{n\binom{n}{k}}\sum_{j=0}^{n-1}\rho^{kj}P(z+\rho^jx).$$
%\end{lemma}
%\textbf{Proof:}
%\begin{align*}
%\sum_{j=0}^{n-1}\rho^{kj}P(z+\rho^jx) &= \sum_{j=0}^{n-1}\rho^{kj}\sum_{l=0}^n\binom{{n}}{l}\check{P}(z^l,x^{n-l})\rho^{(n-l)j}  = \sum_{l=0}^n\binom{{n}}{l}\check{P}(z^l,x^{n-l})\sum_{j=0}^{n-1}\rho^{(k+n-l)j}\\
%&=n\binom{n}{k}\check{P}(z^k,x^{n-k}),
%\end{align*}
%where the last equality holds because, for $0\le l\le n$, we have that  $\sum_{j=0}^{n-1}\rho^{(k+n-l)j}$ is 0 if $l\ne k$ and is equal to $n$ if $l=k$. \hfill$\Box$\medskip

%\begin{theorem}\label{sheet3} Let $E=\ell_p$, $1<p<\infty$, and let $U=B_{\ell_p}$. Then in each sheet there is some $\varphi\in M_b(U)\cap\pi^{-1}(0)$ such that
%%Let $1<p<\infty$ , and let $\varphi\in M_b(B_{\ell_p})$ be such that $\varphi\in\pi^{-1}(0)$.  Then
%$$
%(1-R(\varphi)^p)^{\frac{1}{p}}\cdot B_{\ell_p}\subset \pi(S_{B_{\ell_p}}(\varphi))\subset \big(1-\sup_{m\geq p}\|\varphi_m\|\big)^{1/\lceil p\rceil}\cdot B_{\ell_p},
%$$
% where $\lceil p\rceil$ denotes the smallest natural number which is $\ge p$.
%\end{theorem}
\begin{theorem}\label{sheet3} Let $E=\ell_p$, $1<p<\infty$, and let $U=B_{\ell_p}$. Given a sheet $S$, we take $\varphi\in S\cap\pi^{-1}(0)$ (which exists thanks to Theorem~\ref{sheet2}). Then,
	$$
	(1-R(\varphi)^p)^{\frac{1}{p}}\cdot B_{\ell_p}\subset \pi(S)\subset \big(1-\sup_{m\geq p}\|\varphi_m\|\big)^{1/\lceil p\rceil}\cdot B_{\ell_p},
	$$
	where $\lceil p\rceil$ denotes the smallest natural number which is $\ge p$.
\end{theorem}
\textbf{Proof:} Let us first prove the lower inclusion. Take $z\in(1-R(\varphi)^p)^{\frac{1}{p}}\cdot B_{\ell_p}$. Since $M_b(B_{\ell_p})$ embeds in $M_b(\ell_p)$, we know that $(\varphi_|)^z\in M_b({\ell_p})$. We must show that $(\varphi_|)^z$ belongs to $M_b(B_{\ell_p})$, that is, that $(\varphi_|)^z$ is continuous with respect to the topology in $H_b(B_{\ell_p})$ of uniform convergence on $B_{\ell_p}$-bounded sets. Recall that the seminorms $q_s(f)=\sum_{n=0}^\infty s^n\big\|\frac{d^nf(0)}{n!}\big\|$, with $0<s<1$, define the topology on $H_b(B_{\ell_p})$ (see \cite{din}).

Let $f\in H_b({\ell_p})$ and let us denote by $\sum_{n=0}^{\infty}P_n$ its Taylor series at the origin, then
\begin{align*}
(\varphi_|)^z(f) & =\sum_{n=0}^{\infty}\varphi\big(x\mapsto P_n(x+z)\big).  \end{align*}
Now, since $\|z\|^p+R(\varphi)^p<1$,  we can find $N\in\N$ and $r<1$ such that for every $y\in R(\varphi)\cdot B_{\ell_p}^{(N)}$, we have $z+y\in rB_{\ell_p}$. Then,  by the definition of $R(\varphi)$ and Lemma \ref{lema1}, it follows that
$$
|\varphi\big(x\mapsto P_n(z+x)\big)|\leq \sup_{y\in R(\varphi)\cdot B_{\ell_p}^{(N)}}\|P_n(z+y)\|\leq r^n\|P_n\|.
$$
Therefore,
$$
|(\varphi_|)^z(f)|\le \sum_{n=0}^{\infty}|\varphi\big(x\mapsto P_n(z+x)\big)| \leq \sum_{n=0}^{\infty}r^n\|P_n\|=  q_{r}(f).
$$
 This implies that $(\varphi_|)^z$ belongs to	 $M_b(B_{\ell_p})$.

%
%Let $f\in H_b({\ell_p})$ and let us denote by $\sum_{n=0}^{\infty}P_n$ its Taylor series at the origin, then by Lemma \ref{polarization}
%%using the polarization formula (see \cite[Lemma 1.7]{CarDimMur09} for a proof for $k=1$)
%%$$\check{P}_n(z^k,x^{n-k})=\dfrac{1}{n\binom{n}{k}}\sum_{j=0}^{n-1}\rho^{kj}P_n(z+\rho^jx),\textrm{ where }\rho^n=1\textrm{ and }\rho^j\neq 1\textrm{ for }j<n,$$
%we have
%\begin{align*}
%(\varphi_|)^z(f) & =\sum_{n=0}^{\infty}\varphi\big(x\mapsto P_n(x+z)\big)  =\sum_{n=0}^{\infty}\sum_{k=0}^{n}\binom{n}{k}\varphi\big(x\mapsto\check{P}_n(z^k,x^{n-k})\big) \\
%& =\sum_{n=0}^{\infty}\sum_{k=0}^{n}\sum_{j=0}^{n-1}\rho^{kj}\dfrac{1}{n}\varphi\big(x\mapsto P_n(z+\rho^jx)\big).
%\end{align*}
%
%
%Now, since $\|z\|^p+R(\varphi)^p<1$,  we can find $N\in\N$ and $r<1$ such that for every $x\in R(\varphi)\cdot B_{\ell_p}^{(N)}$, we have $z+x\in rB_{\ell_p}$. Then, by the definition of $R(\varphi)$ and Lemma \ref{lema1}, it follows that
%$$
%|\varphi\big(x\mapsto P_n(z+\rho^jx)\big)|\leq \sup_{y\in R(\varphi)\cdot B_{\ell_p}^{(N)}}\|P_n(z+\rho^jy)\|\leq r^n\|P_n\|.
%$$
%Therefore
%$$
%|(\varphi_|)^z(f)|\le \sum_{n=0}^{\infty}\sum_{k=0}^{n}\sum_{j=0}^{n-1}\dfrac{1}{n}|\varphi\big(x\mapsto P_n(z+\rho^jx)\big)| \leq \sum_{n=0}^{\infty}({n+1})r^n\|P_n\|\leq  cq_{r+\varepsilon}(f),
%$$
%for $\varepsilon>0$ such that $r+\varepsilon<1$ and some constant $c>0$. This implies that $(\varphi_|)^z$ is in $M_b(B_{\ell_p})$.

Now we prove the upper inclusion.
By Theorem \ref{sheet2} we already know that $S_{B_{\ell_p}}(\varphi)$ is a ball centered at zero. Let $z=te_1$, with $t^{\lceil p\rceil}+\sup_{m\geq p}\|\varphi_m\|>1+\delta$, for some $\delta>0$. We will show that  $(\varphi_|)^z$  is not continuous on $H_b(B_{\ell_p})$. This will prove that the radius of the ball $S_{B_{\ell_p}}(\varphi)$ is smaller than or equal to  $(1-\sup_{m\geq p}\|\varphi_m\|)^{1/\lceil p\rceil}$.

Let $0<r<1$ be such that $\varphi\prec rB_{\ell_p}$. Consider $m_0\geq p$ with $t^{\lceil p\rceil}+\|\varphi_{m_0}\|>1+\delta$. For $\varepsilon<{\delta}$, let $P_0\in P(^{m_0}E)$ be such that $\varphi(P_0)>\|\varphi_{m_0}\|-\varepsilon$, and $\|P_0\|\le1$. Note that by Lemma \ref{lema1}, we have that $\varphi(P_0)=\varphi(P_0\circ(I-e_1'\otimes e_1))$. Let $Q_0=P_0\circ(I-e_1'\otimes e_1)$.
%Since $P_0$ is uniformly continuous on $B_{\ell_p}$, there exists $\rho>0$ such that if $x,y\in B_{\ell_p}$ and $\|x-y\|<\rho$ then $|P_0(x)-P_0(y)|<\varepsilon$.
It follows from Lemma \ref{lema1} that
$$
(\varphi_|)^{te_1}(Q_0)=\varphi(x\mapsto Q_0(x+te_1))=\varphi(x\mapsto Q_0(x))=\varphi(Q_0).
$$
%Now $|\varphi^z(P_0)-\varphi(P_0)|=|\varphi^{z^{(N)}}(P_0)-\varphi(P_0)|=|\varphi(P_0(\cdot+z^{(N)})-P_0)|\leq\sup_{x\in rB_{\ell_p}}|P_0(x+z^{(N)})-P_0(x)|<\varepsilon,$ since $\|z^{(N)}\|\leq 1-r$ and then $x+z^{(N)}\in B_{\ell_p}$, and also $\|x+z^{(N)}-x\|<\rho$. Therefore we have $Re(\varphi^{z^{(N)}}(P_0))>\varphi(P_0)-\varepsilon$.
%Let $\theta_k$ be such that $|\theta_k|=1$, for $k=1,\cdots,N$ and $\sum_{k=1}^{N}\theta_k z_k^{\lceil p\rceil}=\sum_{k=1}^{N}|z_k|^{\lceil p\rceil}$, and

Consider the polynomial $Q(x)=(e_1')^{\lceil p\rceil}+Q_0(x)$. Since $m_0\ge p$, we have $\sup_{\|x\|_p\le1}|Q(x)|\leq 1$. Indeed, for $\|x\|_p\le1$,
$$
|Q(x)|\le |x_1|^{\lceil p\rceil}+|P_0\circ(I-e_1'\otimes e_1)(x)|\le |x_1|^{\lceil p\rceil}+\|(I-e_1'\otimes e_1)(x)\|_p^{m_0}\le\|x\|_p^p\le1.
$$
Moreover,
$$
(\varphi_|)^{te_1}(Q)=(\varphi_|)^{te_1}\Big((e_1')^{\lceil p\rceil}+Q_0\Big)=t^{\lceil p\rceil}+\varphi(Q_0),\textrm{ and then}
$$
\begin{align*}
|(\varphi_|)^{te_1}(Q)| = & t^{\lceil p\rceil}+\varphi(Q_0)
>t^{\lceil p\rceil}+\|\varphi_{m_0}\|-\varepsilon>1+\delta-\varepsilon>1+s,\quad\textrm{ for some }s>0.
\end{align*}
Therefore it follows that $|(\varphi_|)^{te_1}(Q^n)|=|(\varphi_|)^{te_1}(Q)|^n>(1+s)^n\tende\infty$ when $n\tende \infty$, while $\|Q^n\|_{B_{\ell_p}}\le1$ for every $n$. Then $(\varphi_|)^{te_1}\notin M_b(B_{\ell_p})$. \hfill$\Box$\medskip

%
%\begin{theorem}\label{sheet4} Let $1<p<\infty$, and let $\varphi\in M_b(B_{\ell_p})$ be such that $\varphi\in\pi^{-1}(0)$. Then $$\pi(S_{B_{\ell_p}}(\varphi))\subset B\Big(0,\big(1-\sup_{m\geq p}\|\varphi_m\|\big)^{1/\lceil p\rceil}\Big),$$ where $\lceil p\rceil$ denotes the smallest natural number which is $\ge p$.\end{theorem}
%
%\textbf{Proof:}  \hfill$\Box$\medskip

%\begin{corollary}\label{sheet5} Let $1<p<\infty$ be a natural number, and let $\varphi\in M_b(B_{\ell_p})$ be such that $\varphi\in\pi^{-1}(0)$. Then $$B(0,(1-R(\varphi)^p)^{\frac{1}{p}})\subset \pi(S_{B_{\ell_p}}(\varphi))\subset B\Big(0,\big(1-\sup_{m\geq p}\|\varphi_m\|\big)^p\Big).$$\end{corollary}

The only homomorphism $\varphi$ such that $\varphi_m=0$ for sufficiently large $m$ is $\delta_0$, so
the previous Theorem allows us to conclude the following.
\begin{corollary} Let $1<p<\infty$, and let $S\subset  M_b(B_{\ell_p})$ be a sheet. Then $\pi(S)=B_{\ell_p}$ if, and only if, $S$ is the sheet of evaluations.\end{corollary}

\begin{remark}\label{remark-summary}\rm
The results of this section can be summarized in the following way. Given a connected component  $S$ of $M_b(B_{\ell_p})$,  there exists $\varphi\in M_b(B_{\ell_p})\cap\pi^{-1}(0)$ and $0<r \le 1$ such that
$$
S= \{\varphi^z\ : \ \|z\|<r\}.
$$
Moreover, $r$ and $\varphi$ satisfy
$$
	(1-R(\varphi)^p)^{\frac{1}{p}} \le  r \le \big(1-\sup_{m\geq p}\|\varphi_m\|\big)^{1/\lceil p\rceil}.
$$

\end{remark}

Some comments deserve to be highlighted. If $p$ is a natural number and $\varphi$ is a homomorphism such that $R(\varphi)=\sup_{m\in\N}\|\varphi_m\|^{\frac{1}{m}}$ is attained at $m=p$ , then it follows that $\pi(S_{B_{\ell_p}}(\varphi))=B(0,(1-R(\varphi)^{p})^{\frac{1}{p}})$, and then we have an accurate description of the sheet of $\varphi$. It is interesting to mention that this is not an artificial hypothesis, since the $r$-block homomorphisms considered in \cite[Definition 5.3]{carmur} satisfy this condition. From this point of view, \cite[Proposition 5.4]{carmur} can be seen now as a consequence of Theorem \ref{sheet3}.

%However, we could not find homomorphisms that are not of $r$-blocks type, and fulfilling the condition $R(\varphi)^p=\|\varphi_p\|$. Nor, we could not find a better description of $\pi(S(\varphi))$ than that made in Corollary \ref{sheet5}. These questions also end up being natural future directions for continuation of this research.

In \cite[Section 6.3]{din}, the spectrum $M_b(E)$ of a symmetrically regular Banach space was informally referred to as the envelope of ``bounded'' holomorphy of $E$ because each bounded type entire function is proved to extend to a holomorphic function on $M_b(E)$ which is of bounded type on each connected component of $M_b(E)$. However, as shown in \cite[Proposition 5.1]{carmur}, the extension need not be of bounded type on the whole Riemann domain, even for a homogeneous polynomial. In the case of the unit ball, we do not know whether the extensions to the spectrum are of bounded type or not. If for any  $\varphi\in M_b(B_{\ell_p})$ the connected components would satisfy $$\pi(S_{B_{\ell_p}}(\varphi))=B(0,(1-R(\varphi)^{p})^{\frac{1}{p}})$$ (that is, if the left inclusion in Theorem~\ref{sheet3} were always an equality), then it would be possible to answer this question affirmatively.

\medskip

\textbf{Acknowledgments}
We would like to thank the anonymous reviewer for his/her comments, which improved the presentation of the article. We are also indebted to Professor Harold P. Boas for carefully reading  the manuscript, for his valuable comments and, in particular, for providing us with a simple proof for Lemma \ref{reinhardt fundamental}.

D. M. Vieira thanks the Departamento de Matem\'{a}tica of the Universidad de Buenos Aires and its members for their kind hospitality.

\end{document}